\documentclass[12pt]{article}
\usepackage{amsmath,amsfonts,amssymb,rotating,colordvi}
\usepackage{colordvi}
\usepackage{mathrsfs}
\def\qed{\nopagebreak\hfill{\rule{4pt}{7pt}}}
\def\proof{\noindent {\it{Proof.} \hskip 2pt}}

\newtheorem{thm}{Theorem}[section]

\newtheorem{prop}[thm]{Proposition}

\newtheorem{conj}[thm]{Conjecture}

\makeatletter \@addtoreset{equation}{section}
\makeatletter
\def\Eqlfill@{\arrowfill@\Relbar\Relbar\Relbar}
\newcommand{\extendEql}[1][]{\ext@arrow 0359\Eqlfill@{#1}}
\makeatother

\numberwithin{equation}{section}

\newdimen\Squaresize \Squaresize=11pt
\newdimen\Thickness \Thickness=0.7pt
\def\Square#1{\hbox{\vrule width \Thickness
   \vbox to \Squaresize{\hrule height \Thickness\vss
    \hbox to \Squaresize{\hss#1\hss}
   \vss\hrule height\Thickness}
\unskip\vrule width \Thickness} \kern-\Thickness}

\def\Vsquare#1{\vbox{\Square{$#1$}}\kern-\Thickness}

\def\moins{\raise 1pt\hbox{{$\scriptstyle -$}}}
\begin{document}

\begin{center}
{\large \bf  The $q$-log-convexity of Domb's polynomials}
\end{center}

\begin{center}
Donna Q. J. Dou$^{1}$ and Anne X. Y. Ren$^{2}$\\[9pt]

$^1$School of Mathematics, Jilin University \\[5pt]
 Changchun 130012, P. R. China\\[9pt]
$^2$Center for Combinatorics, LPMC-TJKLC, Nankai University\\[5pt]
 Tianjin 300071, P. R. China\\[9pt]
$^{1,2}$Center for Applied Mathematics, Tianjin University\\[5pt]
 Tianjin 300072, P. R. China\\[9pt]

Email: $^{1}${\tt
qjdou@jlu.edu.cn}, $^{2}${\tt renxy@nankai.edu.cn}
\end{center}

\noindent\textbf{Abstract.}
In this paper, we prove the $q$-log-convexity of Domb's polynomials, which was conjectured by
Sun in the study of Ramanujan-Sato type series for powers of $\pi$. As a result, we obtain the log-convexity of Domb's numbers.
Our proof is based on the $q$-log-convexity of Narayana polynomials of type $B$ and a criterion for determining $q$-log-convexity of self-reciprocal polynomials.

\noindent \emph{AMS Classification 2010:} 05A20

\noindent \emph{Keywords:} Domb's number, Domb's polynomial, log-convexity, $q$-log-convexity.

\noindent \emph{Corresponding author:} Anne X. Y. Ren, renxy@nankai.edu.cn

\section{Introduction}

The main objective of this paper is to prove a conjecture of Sun \cite{sun2012} on the
$q$-log-convexity of the polynomials $D_n(q)$, which are given by
\begin{align}
D_n(q)&=\sum\limits_{k=0}^n{n\choose k}^2{2k\choose k}{2(n-k)\choose n-k}q^k\label{conj2}.
\end{align}
We call $D_n(q)$ Domb's polynomials because $\{D_n(1)\}_{n\geq0}$ is the ubiquitous sequence of Domb's numbers.
These polynomials $D_n(q)$ were introduced by Sun \cite{sun2012} in his study of the Ramanujan-Sato type series for powers of $\pi$.

Let us first review some definitions. Recall that a nonnegative sequence $\{a_n\}_{n\geq 0}$ is said to be log-concave if, for any $n\geq 1$,
$$a_n^2\geq a_{n-1}a_{n+1};$$
and is said to be log-convex if, for any $n\geq 1$,
$$a_{n-1}a_{n+1} \geq a_n^2.$$
Many sequences arising in combinatorics, algebra and geometry, turn out to be log-concave or log-convex, see Brenti \cite{bren1994} or Stanley \cite{stan1989}.

For a sequence of polynomials with real coefficients, Stanley introduced the notion of $q$-log-concavity.
A polynomial sequence $\{f_n(q)\}_{n\geq 0}$ is said to be $q$-log-concave if, for any $n\geq 1$, the difference
$$f_n^2(q)-f_{n+1}(q)f_{n-1}(q)$$
has nonnegative coefficients as a polynomial in $q$. The $q$-log-concavity of polynomial sequences has been extensively studied, see Bulter \cite{butl1990}, Krattenthaler \cite{krat1989}, Leroux \cite{lero1990} and Sagan \cite{saga1992}. Similarly, a polynomial sequence $\{f_n(q)\}_{n\geq 0}$ is said to be $q$-log-convex if, for any $n\geq 1$, the difference
$$f_{n+1}(q)f_{n-1}(q)-f_n^2(q)$$
has nonnegative coefficients as a polynomial in $q$. Liu and Wang \cite{liu-wang2007} showed that many classical combinatorial polynomials are $q$-log-convex, see also \cite{chen-tang-wang-yang2010, chen-wang-yang2010, chen-wang-yang2011}. It should be noted that Butler and Flanigan \cite{butfla2007} introduced a different kind  of $q$-log-convexity.

This paper is mainly concerned with the log-convexity of Domb's numbers and the $q$-log-convexity of Domb's polynomials.
Domb's numbers appear to be first discovered by Domb \cite{Domb1960}. These numbers play an important role in the study of many subjects, such as Bessel functions \cite{baiborbrogla2008}, random walks \cite{bornuystrwan2011}, interacting systems \cite{Domb1960}, the third order Ap\'ery-like differential equations \cite{almstrzud2011} and the enumeration of abelian squares \cite{ricsha2009}.
Domb's numbers also play an important role in the series expansions for powers of $\pi$. Many new Ramanujan-Sato type series for $1/\pi$ have been derived by using such numbers. For example, Chan, Chan and
Liu \cite{chan-chan-liu2004} obtained a formula for $1/\pi$ involving Domb's numbers:
\begin{align}\label{eq-CCL}
\sum_{n=0}^\infty \frac{5n+1}{64^n}D_n(1)=\frac{8}{\sqrt{3}\pi}.
\end{align}
Sun \cite{sun2012} found over one hundred conjectured series expansion formulas for powers of $\pi$. Among them one conjecture is on the expansion of $1/\pi^2$ in terms
of Domb's numbers,
\[\sum_{n=0}^\infty \frac{40n^2+26n+5}{(-256)^n}{2n\choose n}^2D_n(1)=\frac{24}{\pi^2},\]
see also \cite{sun2011}.

Considering the importance of Domb's numbers and Domb's polynomials, Sun \cite{sun2012-2} further studied their combinatorial properties, and proposed the following conjectures.

\begin{conj}[{\cite[Conjecture 3.12]{sun2012-2}}]\label{sun-conj31}
Both $\{D_{n+1}(1)/D_{n}(1)\}_{n\geq 0}$ and $\{\sqrt[n]{D_{n}(1)}\}_{n\geq 1}$ are strictly increasing. Moreover,  $\{\sqrt[n+1]{D_{n+1}(1)}/\sqrt[n]{D_{n}(1)}\}_{n\geq 1}$ is strictly decreasing.
\end{conj}

\begin{conj}[{\cite[p. 17]{sun2012}}]
The sequence $\{D_n(q)\}_{n\geq 0}$ is $q$-log-convex.
\end{conj}

Note that the log-convexity of $\{D_n(1)\}_{n\geq 0}$ is implied by the strictly increasing property of $\{D_{n+1}(1)/D_{n}(1)\}_{n\geq 0}$. It is also clear that the $q$-log-convexity of $D_n(q)$ also implies the log-convexity of $\{D_n(1)\}_{n\geq 0}$.
Recently, Wang and Zhu \cite{wang-zhu2013} proved the strictly increasing property of $\{D_{n+1}(1)/D_{n}(1)\}_{n\geq 0}$, as well as that of $\{\sqrt[n]{D_{n}(1)}\}_{n\geq 1}$. Chen, Guo and Wang \cite{chen-guo-wang2013} gave a proof of the strictly decreasing property of $\{\sqrt[n+1]{D_{n+1}(1)}/\sqrt[n]{D_{n}(1)}\}_{n\geq 1}$. In this paper, we shall give a proof of the $q$-log-convexity conjecture of  $\{D_n(q)\}_{n\geq 0}$.

It is easy to see that the coefficients of $D_n(q)$ are symmetric.
Such polynomials are called self-reciprocal. More precisely,
we call a polynomial
$$f(q)=a_0+a_1q+\cdots+a_nq^n$$
a self-reciprocal polynomial of degree $n$ if
$f(q)=q^nf(1/q)$. To prove the $q$-log-convexity of $\{D_n(q)\}_{n\geq 0}$, we shall use a criterion for determining $q$-log-convexity of self-reciprocal polynomials, which was first obtained in \cite{douren2012}. To recall this criterion, suppose that $\{f_n(q)\}_{n\geq 0}$
is a $q$-log-convex sequence, where
\begin{align}
f_n(q)=\sum_{k=0}^n a(n,k)q^k.\label{eq-q-log-convex}
\end{align}
For $n\geq 1$ and $0\leq t\leq 2n$, let $\mathcal{L}_t$ be an operator on triangular array $\{a({n, k})\}_{0\leq k\leq n}$ defined by
\begin{align}\label{eq-recur-1}
\mathcal{L}_t(a(n,k))=&a(n+1,k)a(n-1,t-k)+a(n-1,k)a(n+1,t-k)\nonumber\\[5pt]
                     &\qquad \qquad \qquad -2a(n,k)a(n,t-k), \qquad \mbox{ if } 0\leq k\leq \frac{t}{2}.
\end{align}
For the remainder of this paper, unless explicitly stated otherwise, we assume
that $n,k,t$ are nonnegative integers. Then our criterion to be used is as follows.

\begin{thm}[{\cite[Theorem 2.1]{douren2012}}]\label{thm-criterion}
Given a log-convex sequence $\{u_k\}_{k\geq 0}$ and a $q$-log-convex sequence $\{f_n(q)\}_{n\geq 0}$ as defined in \eqref{eq-q-log-convex}, let $\{g_n(q)\}_{n\geq 0}$ be the polynomial sequence defined by
\begin{align}\label{eq-gn1}
g_n(q)=\sum_{k=0}^n a(n,k)u_kq^k.
\end{align}

Assume that the following two conditions are satisfied:
\begin{itemize}
\item[(C1)] for each $n\geq 0$, the polynomial $g_n(q)$ is a self-reciprocal polynomial of degree $n$; and

\item[(C2)] for given $n\geq 1$ and $0\leq t\leq n$, there exists an index $k'$ associated with $n,t$ such that \begin{align*}
\mathcal{L}_t(a(n,k))\left\{
\begin{array}{ll}
\geq 0, & \mbox{ if } 0\leq k\leq k',\\[5pt]
\leq 0, & \mbox{ if } k'< k\leq \frac{t}{2}.
\end{array}
\right.
\end{align*}
\end{itemize}
Then, the polynomial sequence $\{g_n(q)\}_{n\geq 0}$
is $q$-log-convex.
\end{thm}

To use Theorem \ref{thm-criterion} to prove the $q$-log-convexity of $\{D_n(q)\}_{n\geq 0}$, we shall take
\begin{align}\label{Dseq-ua}
u_k=\binom{2k}{k}, \quad a(n,k)={\binom{n}{k}}^2\binom{2n-2k}{n-k},
\end{align}
 and hence, by \eqref{eq-q-log-convex}, \eqref{eq-gn1} and \eqref{Dseq-ua}, we have
\begin{align}\label{Dseq-fg}
f_n(q)=\sum\limits_{k=0}^n{\binom{n}{k}}^2\binom{2n-2k}{n-k}q^k, \quad g_n(q)=D_n(q).
\end{align}
We need to show that the sequence $\{f_n(q)\}_{n\geq 0}$ is $q$-log-convex, and
the triangular array $\{a(n,k)\}_{0\leq k\leq n}$ satisfies the condition (C2) of Theorem \ref{thm-criterion}. For the former, it suffices to prove the $q$-log-convexity of the sequence $\{V_n(q)\}_{n\geq 0}$ given by
\begin{align}\label{v-pol}
V_n(q)=\sum\limits_{k=0}^n{\binom{n}{k}}^2\binom{2k}{k}q^k,
\end{align}
which shall be done in the next section. The key point for proving the latter is to determine the sign changes of certain polynomial of degree 8 over some intervals, which is also the most difficult part for our approach. We complete this task  by examining the properties of the derivatives of this polynomial. To simplify the computations involved, Maple will be frequently used.

\section{The $q$-log-convexity of $\{V_n(q)\}_{n\geq 0}$}

In this section, we are to prove the $q$-log-convexity of the sequence $\{V_n(q)\}_{n\geq 0}$. Our proof is based on the $q$-log-convexity of $\{W_n(q)\}_{n\geq 0}$ given by
\begin{align}\label{nara-B}
W_n(q)=\sum_{k=0}^n{\binom{n}{k}}^2q^k,
\end{align}
which are known as the Narayana polynomials of type $B$.

\begin{thm}[{\cite[Theorem 1.1]{chen-tang-wang-yang2010}}]  \label{q-log-convex-narayana}
The polynomials $W_n(q)$  form a $q$-log-convex sequence.
\end{thm}

The above result was first  conjectured by Liu and Wang \cite{liu-wang2007}, and then proved by Chen, Tang, Wang and Yang \cite{chen-tang-wang-yang2010} by using the theory of symmetric functions. Zhu \cite{zhu2012} gave a simple proof of Theorem \ref{q-log-convex-narayana} based on certain recurrence relation.

To prove the $q$-log-convexity of $\{V_n(q)\}_{n\geq 0}$, we further need a result of Liu and Wang \cite{liu-wang2007}, which provides a mechanism of generating new $q$-log-convex sequences of polynomials from certain log-convex sequences of positive numbers and $q$-log-convex sequences of polynomials.

For $n\geq 1$ and $0\leq t\leq 2n$, let $\mathcal{\widetilde{L}}_t$ be an operator on triangular array $\{a({n, k})\}_{0\leq k\leq n}$ defined by
\begin{align*}
\mathcal{\widetilde{L}}_t(a(n,k))=\left\{
\begin{array}{ll}
a(n+1,k)a(n-1,t-k)+a(n-1,k)a(n+1,t-k)\\[5pt]
\qquad \quad -2a(n,k)a(n,t-k), \qquad \qquad \qquad \mbox{ if }0\leq k< \frac{t}{2},\\[5pt]
a(n+1,k)a(n-1,k)-a^2(n,k), \quad \mbox{ if } t \mbox{ is even and }  k=\frac{t}{2}.
\end{array}
\right.
\end{align*}
Liu and Wang's criterion to determine the $q$-log-convexity of a polynomial sequence is as follows.

\begin{thm}[{\cite[Theorem 4.8]{liu-wang2007}}]\label{liu-wang thm}
Let $\{u_k\}_{k\geq 0}$ be a log-convex sequence and let $\{f_n(q)\}_{n\geq 0}$ be a $q$-log-convex sequence as defined in \eqref{eq-q-log-convex}. Given $n\geq 1$ and $0\leq t\leq 2n$, if there exists an index $k'$ associated with $n,t$ such that
\begin{align*}
\mathcal{\widetilde{L}}_t(a(n,k))\left\{
\begin{array}{ll}
\geq 0, & \mbox{ if }0\leq k\leq k',\\[5pt]
\leq 0, & \mbox{ if } k'< k\leq \frac{t}{2},
\end{array}
\right.
\end{align*}
then, the polynomial sequence $\{g_n(q)\}_{n\geq 0}$ defined by \eqref{eq-gn1} is $q$-log-convex.
%\begin{align}\label{eq-gn}
%g_n(q)=\sum_{k=0}^n a(n,k)u_kq^k
%\end{align}
%is $q$-log-convex.
\end{thm}

The main result of this section is as follows.

\begin{thm} \label{q-log-convex-vpol}
The sequence $\{V_n(q)\}_{n\geq 0}$ given by \eqref{v-pol} is $q$-log-convex.
\end{thm}

\proof Take
\begin{align}
u_k=\binom{2k}{k}, \quad a(n,k)={\binom{n}{k}}^2
\end{align}
in Theorem \ref{liu-wang thm}, and accordingly we have
$$f_n(q)=W_n(q),\qquad g_n(q)=V_n(q).$$
Note that $\{\binom{2k}{k}\}_{k\geq 0}$ is log-convex. By  Theorem
\ref{q-log-convex-narayana}, the sequence $\{f_n(q)\}_{n\geq 0}$ is $q$-log-convex. In the proof of \cite[Theorem 1.3]{chen-tang-wang-yang2010}, it has been shown that
for given $n\geq 1$ and $0\leq t\leq 2n$, there exists
$k'$ such that
\begin{align*}
\mathcal{\widetilde{L}}_t(a(n,k))\left\{
\begin{array}{ll}
\geq 0, & \mbox{ if }0\leq k\leq k',\\[5pt]
\leq 0, & \mbox{ if } k'< k\leq \frac{t}{2}.
\end{array}
\right.
\end{align*}
By Theorem \ref{liu-wang thm}, we obtain the desired $q$-log-convexity of $\{V_n(q)\}_{n\geq 0}$.
\qed

\section{The $q$-log-convexity of $\{D_n(q)\}_{n\geq 0}$}
\label{$q$-log-convexity of $D_n(q)$}
The aim of this section is to prove the $q$-log-convexity of $\{D_n(q)\}_{n\geq 0}$.

\begin{thm}\label{thm-conj3}
The sequence $\{D_n(q)\}_{n\geq 0}$ is $q$-log-convex.
\end{thm}

As discussed before, we plan to use Theorem \ref{thm-criterion} to prove the above theorem. Since  we have obtained the $q$-log-convexity of $\{f_n(q)\}_{n\geq 0}$ as given by \eqref{Dseq-fg}, it remains to prove the following result.

\begin{thm}\label{Dc1 thm3} Let $\{a(n,k)\}_{0\leq k\leq n}$ be the triangular array defined by \eqref{Dseq-ua}, and let $\mathcal{L}_t(a(n,k))$ be given by \eqref{eq-recur-1}.
 Then, for any $n\geq 1$ and $0\leq t\leq n$, there exists an index $k'$ with respect to $n,t$ such that
\begin{align*}
\mathcal{L}_t(a(n,k))\left\{
\begin{array}{ll}
\geq 0, & \mbox{ if }0\leq k\leq k',\\[5pt]
\leq 0, & \mbox{ if } k'< k\leq \frac{t}{2}.
\end{array}
\right.
\end{align*}
\end{thm}

Let us first make some observations. For $n\geq 1$, $0\leq t\leq n$ and $0\leq k\leq t/2$, we have
\begin{align*}
\mathcal{L}_t(a(n,k))=&\binom{n+1}{k}^2\binom{2n-2k+2}{n-k+1}
\binom{n-1}{t-k}^2\binom{2n-2t+2k-2}{n-t+k-1}\nonumber\\[5pt]
&+\binom{n-1}{k}^2\binom{2n-2k-2}{n-k-1}\binom{n+1}{t-k}^2\binom{2n-2t+2k+2}{n-t+k+1}\nonumber\\[5pt]
&-2\binom{n}{k}^2\binom{2n-2k}{n-k}\binom{n}{t-k}^2\binom{2n-2t+2k}{n-t+k}.
\end{align*}
By factorization, we obtain
\begin{align}\label{Deqn3-2}
\mathcal{L}_t(a(n,k))=&\frac{1}{(n-k+1)^3(n-t+k+1)^3(2n-2k-1)(2n-2t+2k-1)}\nonumber\\[5pt]
                       &\times\frac{1}{n^2}\binom{n}{k}^2\binom{2n-2k}{n-k}
                       \binom{n}{t-k}^2\binom{2n-2t+2k}{n-t+k}\psi^{(n,t)}(k),
\end{align}
where
\begin{align}\label{Deqn3-main}
\psi^{(n,t)}(x)=&(n+1)^2(n-x)^3(n-x+1)^3(2n-2t+2x+1)(2n-2t+2x-1)\nonumber\\[5pt]
       &+(n+1)^2(n-t+x)^3(n-t+x+1)^3(2n-2x-1)(2n-2x+1)\nonumber\\[5pt]
       &-2n^2(n-x+1)^3(n-t+x+1)^3(2n-2x-1)(2n-2t+2x-1).
       \end{align}

Clearly, the sign of $\mathcal{L}_t(a(n,k))$ coincides with that of $\psi^{(n,t)}(k)$ unless $t=n$ and $k=0$. Then we could divide the proof of Theorem \ref{Dc1 thm3} into the following three steps:
\begin{itemize}
\item[(S1)] For $n\geq 1$ and $0\leq t\leq n$, prove that $\mathcal{L}_t(a(n,0))\geq 0$ , see Proposition \ref{Dprop-1};

\item[(S2)] For $n\geq 2$ and $0\leq t\leq n-1$, prove that there exists $k'$ such that
\begin{align*}
\psi^{(n,t)}(k)\left\{
\begin{array}{ll}
\geq 0, & \mbox{ if } 1\leq k\leq k',\\[5pt]
\leq 0, & \mbox{ if } k'< k\leq \frac{t}{2},
\end{array}
\right.
\end{align*}
see Proposition \ref{Dprop-2};

\item[(S3)] For $n\geq 2$ and $t=n$, prove that there exists $k'$ such that
\begin{align*}
\psi^{(n,n)}(k)\left\{
\begin{array}{ll}
\geq 0, & \mbox{ if } 1\leq k\leq k',\\[5pt]
\leq 0, & \mbox{ if } k'< k\leq \frac{n}{2},
\end{array}
\right.
\end{align*}
see Proposition \ref{Dprop-3}.
\end{itemize}

We first give a proof of the nonnegativity of $\mathcal{L}_t(a(n,0))$.

\begin{prop}\label{Dprop-1} For any $n\geq 1$ and $0\leq t\leq n$, we have $\mathcal{L}_t(a(n,0))\geq 0$.
\end{prop}
\proof
For $1\leq n\leq 4$, the nonnegativity of $\mathcal{L}_t(a(n,0))$ can be verified directly as follows:
\allowdisplaybreaks
\begin{align*}
\mathcal{L}_0(a(1,0)) &=4,\, \mathcal{L}_1(a(1,0)) =4,\\[5pt]
\mathcal{L}_0(a(2,0)) &=8,\, \mathcal{L}_1(a(2,0)) =32,\, \mathcal{L}_2(a(2,0)) =24,\\[5pt]
\mathcal{L}_0(a(3,0)) &=40,\, \mathcal{L}_1(a(3,0)) =320,\, \mathcal{L}_2(a(3,0)) =646,\, \mathcal{L}_3(a(3,0)) =152,\\[5pt]
\mathcal{L}_0(a(4,0)) &=280,\, \mathcal{L}_1(a(4,0)) =3808,\\[5pt]
\mathcal{L}_2(a(4,0)) &=14296,\,\mathcal{L}_3(a(4,0)) =7772,\, \mathcal{L}_4(a(4,0)) =860.
\end{align*}

Then we assume for the remainder of the proof that $n\geq 5$.
It is easily seen that the sign of $\mathcal{L}_t(a(n,0))$ coincides with that of
\begin{align*}
\frac{\binom{2n}{n}\binom{n}{t}^2\binom{2n-2t}{n-t}\theta(t)}{n^2(n+1)(2n-1)(n-t+1)^3(2n-2t-1)},
\end{align*}
where
\begin{align}\label{Dfunc-g}
\theta(x)=&(2n-1)(2n+1)x^6-3(2n-1)(2n+1)^2x^5\nonumber\\[5pt]
&+(2n-1)(26n^3+41n^2+21n+3)x^4\nonumber\\[5pt]
&-(2n-1)(24n^4+54n^3+44n^2+14n+1)x^3\nonumber\\[5pt]
&+n(4n+1)(n+1)(4n^3+7n^2+3n-3)x^2\nonumber\\[5pt]
&-n^2(8n^2+12n-5)(n+1)^2x+2n^2(2n-1)(n+1)^3.
\end{align}

Then we prove the nonnegativity of $\mathcal{L}_t(a(n,0))$ in two cases:
\begin{itemize}
\item[(i)] When $t=n\geq 5$, we readily checks that
\begin{align*}
\theta(n)&=-n^2(n+1)(n^3+2n^2-3n+2)< 0,
\end{align*}
thus $\mathcal{L}_n(a(n,0))\geq 0$.

\item[(ii)] When $0\leq t<n$, it suffices to show that $\theta(t)\geq 0$. To this end, we consider the monotonicity of $\theta(x)$, regarded as a function of $x$, over the interval $[0,n-1]$.
By \eqref{Dfunc-g}, the derivative of $\theta(x)$ with respect to $x$ is
\begin{align*}
\theta'(x)=&6(2n-1)(2n+1)x^5-15(2n-1)(2n+1)^2x^4\nonumber\\[5pt]
&+4(2n-1)(26n^3+41n^2+21n+3)x^3\nonumber\\[5pt]
&-3(2n-1)(24n^4+54n^3+44n^2+14n+1)x^2\nonumber\\[5pt]
&+2n(4n+1)(n+1)(4n^3+7n^2+3n-3)x\nonumber\\[5pt]
&-n^2(8n^2+12n-5)(n+1)^2.
\end{align*}
We further need the higher order derivatives of $\theta(x)$:
\begin{align*}
\theta''(x)=&30(2n-1)(2n+1)x^4-60(2n-1)(2n+1)^2x^3\nonumber\\[5pt]
&+12(2n-1)(26n^3+41n^2+21n+3)x^2\nonumber\\[5pt]
&-6(2n-1)(24n^4+54n^3+44n^2+14n+1)x\nonumber\\[5pt]
&+2n(4n+1)(n+1)(4n^3+7n^2+3n-3),
\end{align*}
\begin{align*}
\theta'''(x)=&6(2n-1)\left(20(2n+1)x^3-30(2n+1)^2x^2\right.\\[5pt]
&\qquad \qquad \left.+4(26n^3+41n^2+21n+3)x\right.\\[5pt]
&\qquad \qquad \left.-(24n^4+54n^3+44n^2+14n+1)\right),
\end{align*}
\begin{align*}
\theta^{(4)}(x)=&24(2n-1)\left(15(2n+1)x^2-15(2n+1)^2x\right.\\[5pt]
&\qquad \qquad \quad \left.+26n^3+41n^2+21n+3\right).
\end{align*}

Note that the axis of symmetry of the quadratic function $\theta^{(4)}(x)$ is $x=n+\frac{1}{2}$,
then $\theta^{(4)}(x)$ is decreasing on the interval $[0,n-1]$. Since, for $n\geq 5$,
\begin{align*}
\theta^{(4)}(0)&=24(2n-1)(26n^3+41n^2+21n+3)>0,\nonumber\\[5pt] \theta^{(4)}(n-1)&=-24(2n-1)(4n^3+4n^2-66n-33)<0,
\end{align*}
we  deduce that $\theta'''(x)$ increases first and then decreases on $[0,n-1]$.

Further, noting that
\begin{align*}
\theta'''(0)&=-6(2n-1)(24n^4+54n^3+44n^2+14n+1)<0,\nonumber\\[5pt] \theta'''(n-1)&=6(2n-1)(26n^3+26n^2-126n-63)>0,
\end{align*}
we conclude that $\theta''(x)$ first decreases and then increases on the interval $[0,n-1]$.

Moreover, it is easy to check that
\begin{align*}
\theta''(0)&=2n(4n+1)(n+1)(4n^3+7n^2+3n-3)>0,\nonumber\\[5pt]
\theta''(\frac{n}{2})&=-\frac{(68n^5+144n^4-129n^3-244n^2-24n+24)n}{8}<0,\nonumber\\[5pt]
\mbox{and}\, \theta''(n-1)&=8n^6+24n^5-216n^4-112n^3+648n^2-132>0.
\end{align*}
Hence, $\theta'(x)$ increases first, and decreases later, then increases again on $[0,n-1]$.

Furthermore, we find that
\begin{align*}
\theta'(0)&=-n^2(n+1)^2(8n^2+12n-5)<0,\nonumber\\[5pt]
\theta'(1)&=n^2(8n^2-4n-3)(3n^2-8n+1)>0,\nonumber\\[5pt]
\theta'(n-1)&=-8n^6-24n^5+88n^4+48n^3-200n^2+36<0.
\end{align*}
Then, $\theta(x)$ first increases and then decreases on the interval $[1,n-1]$. (Here we do not consider the interval $[0,n-1]$, since we only care about the values of $\theta(t)$ for integer $t$ and we shall deal with $\theta(0)$ separately.)

Finally, we notice that
\begin{align*}
\theta(0)&=2n^2(2n-1)(n+1)^3>0,\nonumber\\[5pt]
\theta(1)&=2n^3(2n-1)(3n^2-3n-2)>0,\nonumber\\[5pt]
\theta(n-1)&=(3n^2+3n-2)(n^4+2n^3-9n^2+6n+4)>0.
\end{align*}
Thus, $\theta(t)>0$ for any integer $0\leq t\leq n-1$.
\end{itemize}
Combining (i) and (ii), we have the desired result.
\qed

Now let us determine the sign of $\psi^{(n,t)}(k)$ for $n\geq 2$ and $0\leq t\leq n-1$.

\begin{prop}\label{Dprop-2}
Given $n\geq 2$ and $0\leq t\leq n-1$,
there exists $k'$ with respect to $n,t$ such that
\begin{align*}
\psi^{(n,t)}(k)\left\{
\begin{array}{ll}
\geq 0, & \mbox{ if } 1\leq k\leq k',\\[5pt]
\leq 0, & \mbox{ if } k'< k\leq \frac{t}{2}.
\end{array}
\right.
\end{align*}
\end{prop}

\proof By \eqref{Deqn3-2} and Proposition \ref{Dprop-1}, we know that
$\psi^{(n,t)}(0)\geq 0$. Therefore, if there exists $t_0: 0\leq t_0\leq t/2$ such that $\psi^{(n,t)}(x)$, regarded as a function of $x$, increases on the interval $[0,t_0)$ and decreases on the interval $[t_0,t/2]$, then Proposition \ref{Dprop-2} should be true. To this end, we need to determine the sign changes of the derivative of $\psi^{(n,t)}(x)$ with respect to $x$ on the interval $[0,t/2]$.

The derivative of $\psi^{(n,t)}(x)$ with respect to $x$ is
\begin{align*}
(\psi^{(n,t)}(x))'&=(2x-t)\psi^{(n,t)}_1(x),
\end{align*}
where
\allowdisplaybreaks
\begin{align*}
\psi^{(n,t)}_1(x)=&32(2n+1)x^6-96(2n+1)tx^5\\[5pt]
&+6x^4\left(32n^4-8n^3(4t-11)+4n^2(2t-7)(t-4)\right.\\[5pt]
&\qquad \quad \left.+2n(24t^2-26t+29)+24t^2-18t+11\right)\\[5pt]
&-4x^3t\left(96n^4-24n^3(4t-11)+12n^2(2t-7)(t-4)\right.\\[5pt]
&\qquad \quad \left.+2n(32t^2-78t+87)+32t^2-54t+33\right)\\[5pt]
&-2x^2\left(128n^6-16n^5(16t-25)+4n^4(12t^2-170t+125)\right.\\[5pt]
&\qquad \quad \left.+2n^3(2t-1)(20t^2+16t-167)\right.\\[5pt]
&\qquad \quad \left.-n^2(28t^4-160t^3+159t^2+341t-134)\right.\\[5pt]
&\qquad \quad \left.-2n(28t^4-82t^3+72t^2+38t-17)\right.\\[5pt]
&\qquad \quad \left.-28t^4+66t^3-36t^2-11t+6\right)\\[5pt]
&+2xt(n+1)\left(128n^5-16n^4(16t-17)+4n^3(36t^2-106t+57)\right.\\[5pt]
&\qquad \quad \qquad \quad \left.-2n^2(8t^3-72t^2+138t-53)\right.\\[5pt]
&\qquad \quad \qquad \quad \left.-n(4t^4+4t^3-33t^2+65t-28)\right.\\[5pt]
&\qquad \quad \qquad \quad \left.-4t^4+12t^3-3t^2-11t+6\right)\\[5pt]
&+(n+1)^2\left(64n^6-16n^5(12t-5)+8n^4(22t^2-21t-3)\right.\\[5pt]
&\qquad \quad \qquad \left.-8n^3(4t^3-6t^2-9t+7)\right.\\[5pt]
&\qquad \quad \qquad \left.-2n^2(12t^4-32t^3+51t^2-45t+11)\right.\\[5pt]
&\qquad \quad \qquad \left.+2n(t-1)(4t^4-8t^3+19t^2-15t+3)\right.\\[5pt]
&\qquad \quad \qquad \left.-6t^4+15t^3-12t^2+3t\right).
\end{align*}
Furthermore, the derivative of $\psi^{(n,t)}_1(x)$ is:
\begin{align}\label{Deq-diffh1}
(\psi^{(n,t)}_1(x))'&=2(2x-t)\psi^{(n,t)}_2(x),
\end{align}
where
\allowdisplaybreaks
\begin{align*}
\psi^{(n,t)}_2(x)=&48x^4(2n+1)-96x^3t(2n+1)\\[5pt]
&+6x^2\left(32n^4-8n^3(4t-11)+4n^2(2t-7)(t-4)\right.\\[5pt]
&\qquad \quad \left.+2n(16t^2-26t+29)+16t^2-18t+11\right)\\[5pt]
&-6xt\left(32n^4-8n^3(4t-11)+4n^2(2t-7)(t-4)\right.\\[5pt]
&\qquad \quad \left.+2n(8t^2-26t+29)+8t^2-18t+11\right)\\[5pt]
&-(n+1)\left(128n^5-16n^4(16t-17)+4n^3(36t^2-106t+57)\right.\\[5pt]
&\qquad \quad \qquad \left.-2n^2(8t^3-72t^2+138t-53)\right.\\[5pt]
&\qquad \quad \qquad \left.-(4t^4+4t^3-33t^2+65t-28)n\right.\\[5pt]
&\qquad \quad \qquad \left.-4t^4+12t^3-3t^2-11t+6\right).
\end{align*}
The derivative of $\psi^{(n,t)}_2(x)$ is:
\begin{align}
(\psi^{(n,t)}_2(x))'&=6(2x-t)\psi^{(n,t)}_3(x),
\end{align}
where
\allowdisplaybreaks
\begin{align*}
\psi^{(n,t)}_3(x)=&16x^2(2n+1)-16xt(2n+1)\\[5pt]
&+32n^4-8n^3(4t-11)+4n^2(2t-7)(t-4)\\[5pt]
&\qquad \quad +2n(8t^2-26t+29)+8t^2-18t+11.
\end{align*}
Note that the axis of symmetry of the quadratic function  $\psi^{(n,t)}_3(x)$ is $x=t/2$. Hence, $\psi^{(n,t)}_3(x)$ decreases as $x$ increases from $0$ to $t/2$. In addition, for $n\geq 1$ and $0\leq t\leq n$, we have
\begin{align*}
\psi^{(n,t)}_3\left(\frac{t}{2}\right)=&(8n^2+8n+4)(n-t)^2+(16n^3+44n^2+44n+18)(n-t)\\[5pt]
&+8n^4+36n^3+64n^2+40n+11>0,
\end{align*}
then $\psi^{(n,t)}_2(x)$ decreases on the interval $[0, t/2]$.

We also check that, for $n\geq 1$ and $0\leq t<n$,
\begin{align*}
\psi^{(n,t)}_2\left(\frac{t}{2}\right)=&-(8n^2+10n+5)(n-t)^4-(32n^3+70n^2+50n+15)(n-t)^3\\[5pt]
&-(48n^4+150n^3+165n^2+72n+\frac{27}{2})(n-t)^2\\[5pt]
&-(32n^5+130n^4+200n^3+152n^2+49n+11)(n-t)\\[5pt]
&-(8n^6+40n^5+80n^4+95n^3+\frac{143}{2}n^2+23n+6)<0,\\[5pt]
\psi^{(n,t)}_2\left(0\right)=&(n+1)\left(4t^4(n+1)+4t^3(n+1)(4n-3)\right.\\[5pt]
&\qquad \qquad \left.-3t^2(48n^3+48n^2+11n-1)\right.\\[5pt]
&\qquad \qquad \left.+t(256n^4+424n^3+276n^2+65n+11)\right.\\[5pt]
&\qquad \qquad \left.-2(n+1)(64n^4+72n^3+42n^2+11n+3)\right).
\end{align*}
So when $\psi^{(n,t)}_2(0)>0$, there exists $0<x_0<t/2$ satisfying
\begin{align}\label{Deq-psi2t0}
\psi^{(n,t)}_2(x_0)=0,
\end{align}
and $\psi^{(n,t)}_1(x)$ is decreasing on $[0,x_0]$ and increasing on $[x_0,t/2]$.

Otherwise, $\psi^{(n,t)}_2(0)\leq0$, $\psi^{(n,t)}_1(x)$ is increasing on $[0,t/2]$.

We are now ready to determine the sign changes of $(\psi^{(n,t)}(x))'$ based on the above monotonicity of $\psi^{(n,t)}_1(x)$. For this purpose, the values of $\psi^{(n,t)}_1(x)$ at the two endpoints of the interval $[0,t/2]$ are to be examined.

It is easily verified that $\psi^{(n,t)}_1(t/2)>0$ for any integers $n\geq 2$ and $0\leq t< n$. Using Maple, we find that
\begin{align*}
\psi^{(n,t)}_1\left(\frac{t}{2}\right) =&\left(4(2n-t)^5(2n^2+2n+1)\right.\\[5pt]
&\left.+2(2n-t)^4(10n^2-2n-1)\right.\\[5pt]
&\left.+(2n-t)^3(20n^2-46n-23)\right.\\[5pt]
&\left.+10(2n-t)^2(2n^2-6n-3)\right.\\[5pt]
&\left.+4(2n-t)(14n^2-6n-3)\right.\\[5pt]
&\left.+56n^2\right)\frac{(2n-t+2)}{8},
\end{align*}
which is greater than $0$ whenever $n\geq 4$ and $0\leq t<n$. It remains to check that $\psi^{(n,t)}_1(t/2)>0$ for $n=2,3$. In
fact, for $n=2$, we have $0\leq t<2$ and hence
\begin{align*}
\psi^{(2,t)}_1\left(\frac{t}{2}\right)=&(52(4-t)^5+70(4-t)^4-35(4-t)^3-70(4-t)^2+880-164t)\\[5pt]
&\times\frac{6-t}{8}>0.
\end{align*}
For $n=3$, we have $0\leq t<3$ and hence
\begin{align*}
\psi^{(3,t)}_1\left(\frac{t}{2}\right)=&(100(6-t)^5+166(6-t)^4+19(6-t)^3\\[5pt]
&-30(6-t)^2-420t+3024))\\[5pt]
&\times\frac{8-t}{8}>0.
\end{align*}
As we see, the value of $\psi^{(n,t)}_1(t/2)$ must be positive. We consider the monotonicity of $\psi^{(n,t)}(x)$ in the following three cases by taking into account the value of $x_0$ and the sign of $\psi^{(n,t)}_1(0)$:
\begin{itemize}
\item[(i)] $\psi^{(n,t)}_1(x)$ increases on $[0,t/2]$ and $\psi^{(n,t)}_1(0)\geq 0$. In this case, $\psi^{(n,t)}_1(x)$ increases from a nonnegative value to a positive value as $x$ increases from $0$ to $t/2$. Thus, $(\psi^{(n,t)}(x))'$ takes only nonpositive values on $[0,t/2]$. This means that $\psi^{(n,t)}(x)$ decreases on the interval $[0,t/2]$.

\item[(ii)] $\psi^{(n,t)}_1(x)$ increases on $[0,t/2]$ and $\psi^{(n,t)}_1(0)<0$. In this case, $\psi^{(n,t)}_1(x)$ increases from a negative value to a positive value as $x$ increases from $0$ to $t/2$. Then, there exists $0<t_0<t/2$ satisfying
    \begin{align*}
      \psi^{(n,t)}_1(x)\left\{
      \begin{array}{ll}
         \leq 0,& \mbox{ if } 0\leq x\leq t_0,\\[5pt]
         \geq 0,& \mbox{ if } t_0<x\leq t/2.
      \end{array}
      \right.
    \end{align*}
Hence,   we have
  \begin{align*}
      (\psi^{(n,t)}(x))'\left\{
      \begin{array}{ll}
         \geq 0,& \mbox{ if } 0\leq x\leq t_0,\\[5pt]
         \leq 0,& \mbox{ if } t_0<x\leq t/2,
      \end{array}
      \right.
    \end{align*}
which means that $\psi^{(n,t)}(x)$ is increasing on $[0,t_0]$ and decreasing on $[t_0,t/2]$.

\item[(iii)] $\psi^{(n,t)}_1(x)$ is decreasing on $[0,x_0]$ and increasing on $[x_0,t/2]$ ($x_0$ is defined in \eqref{Deq-psi2t0}). In this case, we must have $\psi^{(n,t)}_1(0)< 0$. Once this assertion is proved, we obtain the desired monotonicity of $\psi^{(n,t)}(x)$ on $[0,t/2]$, by using similar arguments as in case (ii).
    Note that $\psi^{(n,t)}_2(0)>0$ in this case.

    Now we are to deduce $\psi^{(n,t)}_1(0)< 0$ from $\psi^{(n,t)}_2(0)>0$.
    Using Maple, we find that
    \begin{align*}
    \psi^{(n,t)}_1(0)=&(n+1)^2\left(8nt^5-6t^4(2n+1)^2-t^3(32n^3-64n^2-54n-15)\right.\\[5pt]
                      &\qquad \qquad \left.+2t^2(88n^4+24n^3-51n^2-34n-6)\right.\\[5pt]
                      &\qquad \qquad \left.-3t(64n^5+56n^4-24n^3-30n^2-12n-1)\right.\\[5pt]
                      &\qquad \qquad \left.+2n(n+1)(32n^4+8n^3-20n^2-8n-3)\right),
    \end{align*}
\begin{align*}
    \psi^{(n,t)}_2\left(0\right)=&(n+1)\left(4t^4(n+1)+4t^3(n+1)(4n-3)\right.\\[5pt]
                      &\qquad \qquad \left.-3t^2(48n^3+48n^2+11n-1)\right.\\[5pt]
                      &\qquad \qquad \left.+t(256n^4+424n^3+276n^2+65n+11)\right.\\[5pt]
                      &\qquad \qquad \left.-2(n+1)(64n^4+72n^3+42n^2+11n+3)\right).
    \end{align*}
    Recall that $0\leq t\leq n-1$ by the hypothesis.
    We may regard $\psi^{(n,t)}_1(0)/(n+1)^2$ as a polynomial in the variable $t$ over the interval $[0,n-1]$, denoted by $\xi(t)$, and similarly, regard $\psi^{(n,t)}_2(0)/(n+1)$ as a polynomial $\eta(t)$. Now we can divide the proof of $\psi^{(n,t)}_1(0)< 0$ into the following three statements:

    \medskip

    \textbf{Claim 1.} For $0\leq t\leq n-1$, if $\psi^{(n,t)}_2(0)>0$, then $n\neq 2,3,4$.

    \noindent \textit{Proof of Claim 1.} In fact, it is routine to check that $\psi^{(n,t)}_2(0)<0$ if $(n,t)\in \{(2,0),(2,1),(3,0),(3,1),(3,2),(4,0),(4,1),(4,2),(4,3)\}$, which contradicts the positivity of $\psi^{(n,t)}_2(0)$. \qed

    \medskip

    \noindent \textbf{Claim 2.} For any integer $n\geq 4$ and $0\leq t\leq n-1$, the polynomial $\xi(t)$ takes only negative values on the interval $[\frac{3}{4}n,n-1]$.

    \noindent \textit{Proof of Claim 2.}
    Note that, for $4\leq n\leq 7$, $t$ only takes the integer value $n-1$ on the interval $[\frac{3}{4}n,n-1]$.
    And it is routine to check that for $n\geq 4$,
    \begin{align*}
    \xi(n-1)&=-8n^5-14n^4+119n^3+75n^2-100n-36<0.
    \end{align*}
    When $n\geq 8$, we need to consider the derivatives of the first three orders of $\xi(t)$ with respect to $t$:
    \begin{align*}
      \xi'(t)=&40nt^4-24t^3(2n+1)^2-3t^2(32n^3-64n^2-54n-15)\\[5pt]
        &\quad +4t(88n^4+24n^3-51n^2-34n-6)\\[5pt]
        &\quad -192n^5-168n^4+72n^3+90n^2+36n+3,\\[5pt]
      \xi''(t)=&160nt^3-72t^2(2n+1)^2-6t(32n^3-64n^2-54n-15)\\[5pt]
        &\quad +352n^4+96n^3-204n^2-136n-24,\\[5pt]
      \xi'''(t)=&480nt^2-144t(2n+1)^2-192n^3+384n^2+324n+90.
    \end{align*}
    Note that, when $n\geq 8$, the axis of symmetry of the quadratic function $\xi'''(t)$ is $t=\frac{3(2n+1)^2}{20n}<\frac{3n}{4}$, meaning that $\xi'''(t)$ increases on the interval $[\frac{3}{4}n,n-1]$.

    Further
      $$\xi'''(n-1)=-288n^3-576n^2+1236n+234<0,$$
    then $\xi'''(t)<0$ when $t\in[\frac{3}{4}n,n-1]$ and $\xi''(t)$ decreases on the interval $[\frac{3}{4}n,n-1]$.

    Similarly, since
     $$\xi''(n-1)=32n^4+480n^3+432n^2-674n-186>0,$$
    $\xi'(t)$ increases on the interval $[\frac{3}{4}n,n-1]$.

    Also,
    \begin{align*}
     \xi'(\frac{3n}{4})=&-\frac{315}{32}n^5-\frac{57}{2}n^4+\frac{213}{16}n^2+18n+3<0,\\[5pt]
    \mbox{and}\ \xi'(n-1)=&8n^5-8n^4-330n^3-209n^2+284n+96>0,
    \end{align*}
    imply that $\xi(t)$ decreases at first and then increases as $t$ varies from $\frac{3n}{4}$ to $n-1$.

    Recall that $\xi(n-1)<0$ and
    \begin{align*}
     \xi\left(\frac{3}{4}n\right)&=-\frac{n}{128}(25n^5-52n^4+831n^3+2614n^2+224n+480)<0,
    \end{align*}
    we obtain $\xi(t)<0$ for any $t\in[\frac{3}{4}n,n-1]$ when $n\geq8$.

    This completes the proof of Claim 2.\qed

\medskip

\noindent \textbf{Claim 3.} For any integer $n\geq 4$ and $0\leq t\leq n-1$, $\eta(t)<0$ for any  $t\in[0,\frac{3}{4}n]$.

\noindent \textit{Proof of Claim 3.} For $n\geq 4$, it is easy to compute that
\begin{align*}
\eta(0)&=-2(n+1)(64n^4+72n^3+42n^2+11n+3)<0,\\[5pt]
\eta\left(\frac{3}{4}n\right)&=-\frac{575}{64}n^5-\frac{2051}{64}n^4-\frac{357}{8}n^3-\frac{889}{16}n^2-\frac{79}{4}n-6<0.
\end{align*}
The first order derivative and the second order derivative of $\eta(t)$ with respect to $t$ are
\begin{align*}
\eta'(t)=&16t^3(n+1)+12t^2(n+1)(4n-3)-6t(48n^3+48n^2+11n-1)\\[5pt]
&\quad +256n^4+424n^3+276n^2+65n+11,\\[5pt]
\eta''(t)=&48t^2(n+1)+24t(n+1)(4n-3)-288n^3-288n^2-66n+6.
\end{align*}
Note that the axis of symmetry of the quadratic function $\eta''(t)$ is $$t=-(n-\frac{3}{4})<0,$$
that is to say, $\eta''(t)$ increases on the interval $[0,3n/4]$.

And with,
\begin{align*}
\eta''(3n/4)=&-189n^3-243n^2-120n+6<0,
\end{align*}
we have $\eta''(t)<0$ for any $t\in [0,3n/4]$, then $\eta'(t)$ decreases as $t$ increases from 0 to $3n/4$.

From
\begin{align*}
\eta'(3n/4)=&\frac{n+1}{4}(295n^3+591n^2+234n+44)>0,
\end{align*}
$\eta(t)$ should be increasing on the interval $[0,3n/4]$.
Therefore, $\eta(t)<0$ for any $t\in [0,3n/4]$ since both $\eta(0)$ and $\eta(3n/4)$ are negative. This ends the proof of Claim 3.\qed

\medskip
Now we can prove $\psi^{(n,t)}_1(0)<0$. Since $\psi^{(n,t)}_2(0)>0$, it follows that $\eta(t)>0$ and $n\geq 5$ by Claim 1, then by Claim 3, we must have $t> 3n/4$. Then by Claim 2, we get $\xi(t)<0$, and hence $\psi^{(n,t)}_1(0)<0$, as desired.
\end{itemize}
Combining (i), (ii) and (iii), we have completed the proof.\qed

The above proposition is the key step for the proof of  Theorem \ref{Dc1 thm3}. But we need one more proposition.

\begin{prop}\label{Dprop-3}
Given $n\geq 2$, there exists $k'$ with respect to $n$ such that
\begin{align*}
\psi^{(n,n)}(k)\left\{
\begin{array}{ll}
\geq 0, & \mbox{ if } 1\leq k\leq k',\\[5pt]
\leq 0, & \mbox{ if } k'< k\leq \frac{n}{2}.
\end{array}
\right.
\end{align*}
\end{prop}

By \eqref{Deqn3-main}, we have
\begin{align*}
\psi^{(n,n)}(x)=&8x^8(2n+1)-32nx^7(2n+1)\\[5pt]
&+2x^6(8n^4+92n^3+92n^2+40n+11)\\[5pt]
    &-2nx^5(24n^4+164n^3+220n^2+120n+33)\\[5pt]
    &+x^4(52n^6+300n^5+435n^4+205n^3+11n^2-23n-6)\\[5pt]
    &-2nx^3(12n^6+64n^5+87n^4+5n^3-44n^2-23n-6)\\[5pt]
    &+nx^2(n+1)(4n^6+16n^5-3n^4-63n^3-34n^2-7n-3)\\[5pt]
    &+n^2x(6n^3+19n^2-2n+3)(n+1)^2\\[5pt]
    &-n^2(n^3+2n^2-3n+2)(n+1)^3.
\end{align*}
To determine the sign of $\psi^{(n,n)}(k)$, the derivative of $\psi^{(n,n)}(x)$ with respect to $x$ would be considered. Using Maple, we have
\begin{align*}
(\psi^{(n,n)}(x))'&=(2x-n)\psi^{(n,n)}_1(x),
\end{align*}
where
\begin{align*}
\psi^{(n,n)}_1(x)=&32x^6(2n+1)-96nx^5(2n+1)\\[5pt]
&+6x^4(8n^4+76n^3+84n^2+40n+11)\\[5pt]
&-4nx^3(24n^4+148n^3+212n^2+120n+33)\\[5pt]
&+2x^2(28n^6+152n^5+223n^4+85n^3-22n^2-23n-6)\\[5pt]
&-2nx(n+1)(4n^5+16n^4+3n^3-38n^2-17n-6)\\[5pt]
&-n(6n^3+19n^2-2n+3)(n+1)^2.
\end{align*}
Taking the derivative of $\psi^{(n,n)}_1(x)$ with respect to $x$ again:
\begin{align*}
(\psi^{(n,n)}_1(x))'&=2(2x-n)\psi^{(n,n)}_2(x),
\end{align*}
where
\begin{align*}
\psi^{(n,n)}_2(x)=&48x^4(2n+1)-96nx^3(2n+1)\\[5pt]
&+6x^2(8n^4+60n^3+76n^2+40n+11)\\[5pt]
&-6nx(8n^4+44n^3+68n^2+40n+11)\\[5pt]
&+(n+1)(4n^5+16n^4+3n^3-38n^2-17n-6).
\end{align*}

We further need to consider the derivative of $\psi^{(n,n)}_2(x)$ with respect to $x$:
\begin{align*}
(\psi^{(n,n)}_2(x))'&=6(2x-n)\psi^{(n,n)}_3(x),
\end{align*}
where
\begin{align*}
\psi^{(n,n)}_3(x)=&16x^2(2n+1)-16nx(2n+1)+8n^4+44n^3+68n^2+40n+11.
\end{align*}
Note that the axis of symmetry of the quadratic function $\psi^{(n,n)}_3(x)$ is $x=n/2$, and, for $n\geq 2$,
\begin{align*}
\psi^{(n,n)}_3(0)&=8n^4+44n^3+68n^2+40n+11>0,\\[5pt]
\psi^{(n,n)}_3(n/2)&=8n^4+36n^3+64n^2+40n+11>0.
\end{align*}
Thus, $\psi^{(n,n)}_3(x)$ decreases from a positive value to a positive value as $x$ increases from $0$ to ${n}/{2}$,
which implies that $\psi^{(n,n)}_2(x)$ decreases as $x$ increases from $0$ to ${n}/{2}$.

In addition, for $n\geq 2$,
\begin{align*}
\psi^{(n,n)}_2(0)&=4n^6+20n^5+19n^4-35n^3-55n^2-23n-6>0,\\[5pt]
\psi^{(n,n)}_2(n/2)&=-8n^6-40n^5-80n^4-95n^3-\frac{143}{2}n^2-23n-6<0,
\end{align*}
then $\psi^{(n,n)}_2(x)$ decreases from a positive value to a negative value as $x$ increases from $0$ to ${n}/{2}$.
Hence, there exists $0<x_0<n/2$ such that
    \begin{align*}
      (\psi^{(n,n)}_1(x))'\left\{
      \begin{array}{ll}
         \leq 0,& \mbox{ if } 0\leq x\leq x_0,\\[5pt]
         \geq 0,& \mbox{ if } x_0<x\leq n/2.
      \end{array}
      \right.
    \end{align*}
In view of that, for $n\geq 2$,
\begin{align*}
\psi^{(n,n)}_1(0)&=-6n^6-31n^5-42n^4-18n^3-4n^2-3n<0,\\[5pt]
\psi^{(n,n)}_1(n/2)&=n^8+\frac{11}{2}n^7+\frac{19}{2}n^6+\frac{3}{2}n^5-\frac{83}{8}n^4-\frac{13}{2}n^3-n^2-3n>0,
\end{align*}
there exists $0< x_1<n/2$ with
    \begin{align*}
      \psi^{(n,n)}_1(x)\left\{
      \begin{array}{ll}
         \leq 0,& \mbox{ if } 0\leq x\leq x_1,\\[5pt]
         \geq 0,& \mbox{ if } x_1<x\leq n/2.
      \end{array}
      \right.
    \end{align*}
Therefore,
    \begin{align*}
      (\psi^{(n,n)}(x))'\left\{
      \begin{array}{ll}
         \geq 0,& \mbox{ if } 0\leq x\leq x_1,\\[5pt]
         \leq 0,& \mbox{ if } x_1<x\leq n/2.
      \end{array}
      \right.
    \end{align*}
Moreover, it is easy to verify that, for $n\geq 2$,
\begin{align*}
\psi^{(n,n)}(0)&=-n^2(n^3+2n^2-3n+2)(n+1)^3<0, \\[5pt]
\psi^{(n,n)}(1)&=n^6(3n^2-3n-38)+3n^2(18n^3-n-24)+35n^4-16n+24>0, \\[5pt]
\psi^{(n,n)}(n/2)&=-\frac{1}{32}n^2(n-1)(2n^3+3n^2-5n-8)(n+2)^3<0.
\end{align*}
Then, there exists $1< x_2<n/2$ satifying
    \begin{align*}
      \psi^{(n,n)}(x)\left\{
      \begin{array}{ll}
         \geq 0,& \mbox{ if } 1\leq x\leq x_2,\\[5pt]
         \leq 0,& \mbox{ if } x_2<x\leq n/2.
      \end{array}
      \right.
    \end{align*}
Thus, there exists an index $k'=k'(n,n)$ for which $\psi^{(n,n)}(k)\geq 0$ for $1\leq k\leq k'$ and $\psi^{(n,n)}(k)\leq 0$ for $k'<k\leq n/2$, as desired. This completes the proof.\qed

We now come to the proof of Theorem \ref{Dc1 thm3}.

\noindent \textit{Proof of Theorem \ref{Dc1 thm3}.} By Proposition \ref{Dprop-1}, for any $n\geq 1$ and $0\leq t\leq n$, we have $\mathcal{L}_t(a(n,0))\geq 0$.  Given $n\geq 1$, it suffices to show that, for $0\leq t\leq n$, there exists $k'$ such that $\mathcal{L}_t(a(n,k))\geq 0$ for $1\leq k\leq k'$ and $\mathcal{L}_t(a(n,k))\leq 0$ for $k'< k\leq t/2$. By \eqref{Deqn3-2}, for $k\geq 1$, the sign of $\mathcal{L}_t(a(n,k))$ coincides with that of $\psi^{(n,t)}(k)$. From Propositions \ref{Dprop-2} and \ref{Dprop-3}, we obtain the desired result. \qed

Finally, we can prove the $q$-log-convexity of $\{D_n(q)\}_{n\geq 0}$.

\noindent \textit{Proof of Theorem \ref{thm-conj3}.} Combining Theorems \ref{thm-criterion}, \ref{q-log-convex-vpol} and \ref{Dc1 thm3}, we obtain the desired result. \qed

\vskip 3mm
\noindent {\bf Acknowledgments.} This work was supported by the 973 Project, the PCSIRT Project
of the Ministry of Education, the National Science Foundation of China, and the Fundamental Research Funds for the Central Universities of China.

\end{document}